# Harmonic Functions for Data Reconstruction on 3D Manifolds

Li Chen
Department of Computer Science and Information Technology
University of the District of Columbia
Email: *lchen@udc.edu*

Feng Luo
Department of Mathematics
Rutgers University, New Brunswick
Email: *fluo@math.rutgers.edu*

**Abstract:** In computer graphics, smooth data reconstruction on 2D or 3D manifolds usually refers to subdivision problems. Such a method is only valid based on dense sample points. The manifold usually needs to be triangulated into meshes (or patches) and each node on the mesh will have an initial value. While the mesh is refined the algorithm will provide a smooth function on the redefined manifolds.

However, when data points are not dense and the original mesh is not allowed to be changed, how is the "continuous and/or smooth" reconstruction possible? This paper will present a new method using harmonic functions to solve the problem. Our method contains the following steps: (1) Partition the boundary surfaces of the 3D manifold based on sample points so that each sample point is on the edge of the partition. (2) Use gradually varied interpolation on the edges so that each point on edge will be assigned a value. In addition, all values on the edge are gradually varied. (3) Use discrete harmonic function to fit the unknown points, i.e. the points inside each partition patch.

The fitted function will be a harmonic or a local harmonic function in each partitioned area. The function on edge will be "near" continuous (or "near" gradually varied). If we need a smoothed surface on the manifold, we can apply subdivision algorithms.

This paper has also a philosophical advantage over triangulation meshes. People usually use triangulation for data reconstruction. This paper employs harmonic functions, a generalization of triangulation because linearity is a form of harmonic. Therefore, local harmonic initialization is more sophisticated then triangulation. This paper is a conceptual and methodological paper. This paper does not focus on detailed mathematical analysis nor fine algorithm design.

## 1. Introduction

To get a smooth function on a 2D or 3D manifold is a common problem in computer graphics and computational mathematics. In computer graphics, smooth data reconstruction on 2D or 3D manifolds usually refers to subdivision problems. Such a method is only valid based on dense sample points. The manifold usually needs to be triangulated into meshes (or patches) and each node on the mesh will have an initial value. While the mesh is refined the algorithm will provide a smooth function on the redefined manifolds.

When data points are not dense and the original mesh is not supposed to be changed, how to make a "continuous and/or smooth" reconstruction?

Recently, we have developed the algorithms for smooth-continuous data reconstruction for a 2D region based on the so-called digital-discrete method. The method uses gradually varied functions in each smooth order level to obtain a smooth function in a 2D region. We also used gradually varied function on manifold to reconstruct a continuous surface on the manifold.

However, when data points are not dense and the original mesh is not allowed to be changed, how is the "continuous and/or smooth" reconstruction possible? This paper will present a new method using harmonic functions to solve the problem. Our method contains the following steps: (1) Partition the boundary surfaces of the 3D manifold based on sample points so that each sample point is on the edge of the partition. (2) Use gradually varied interpolation on the edges so that each point on edge will be assigned a value. In addition, all values on the edge are gradually varied. (3) Use discrete harmonic function to fit the unknown points, i.e. the points inside each partition patch. Finally, we can use the boundary surface to do the harmonic reconstruction for the original 3D manifold.

The fitted function will be a harmonic or a local harmonic function in each partitioned area. The function on edge will be "near" continuous (or "near" gradually varied). If we need a smoothed surface on the manifold, we can apply subdivision algorithms.

The new method we present may eliminate the use of triangulations, one of the foundations of computer graphics for at least 30 years.

In the past, people usually use triangulation for data reconstruction. This paper employs harmonic functions, a generalization of triangulation, because linearity is a form of harmonic. Therefore, local harmonic initialization is more sophisticated then triangulation.

This paper contains a philosophical change in computer graphics. Triangulation is no longer needed with this new method, a discovery which is incredibly exciting! Moreover, the method we present is the generalization of triangulation and bi-linear interpolation (in Coons surfaces, 4 points interpolation).

It is true that I did not realize before submitting the paper that triangulation and bilinear are both harmonic. However, harmonic functions do not only refer to linear functions. It is much more general.

We represent major steps of the process of the data reconstruction described in this paper as follows:

**Step 1:** Given a set of sample points on a manifold, the objective is to get a data fitting based on the sample points.
**Step 2:** Split the boundary surface of the manifold into surface pieces (polygons or any shape) such that all samples are on the edge of the surfaces pieces.
**Step 3:** Find the gradually varied fitting at each point on the edge of all surface pieces (polygons). This resolution will guarantee the continuity of global fitting on the manifold. The traditional existing method cannot make this possible.
**Step 4:** Make harmonic function fitting in each surface piece, inside of each polygon.

## 2. Discrete Harmonic Functions and Gradually Varied Functions

In this section, we introduce some concepts used in this paper. Let G=(D, E) be a graph.

Discrete Harmonic Functions

Let $f: D \rightarrow R$ be a function on $D$. $f$ is said to be discretely harmonic if $f(x)$ is equal to the averge value for all adjacent vertices of x.

**Theorem 2.1** (Harmonic Functions) For a bounded region $D$ and its boundary J, if $f$ on $J$ is continuous, there is a unique harmonic extension $F$ of $f$ such that the extension is harmonic in $D$-$J$.

More generally, the solution for the Dirichlet Problem:

**Theorem 2.1'** [1] For a 2D piece-wise linear manifold with multiple closed simple curves as it's boundary . (The genus of the manifold could be $g>0$.) If the function $f$ on each of these curves is continuous, then there will be a harmonic extension of the function $f$ on the whole manifold.

This theorem is valid in discrete domain in terms of approximation. We need to give a reasonable interpretation of continuous function in discrete case. In this paper, we prefer to use gradually varied functions.

Gradually Varied Function (GVF) Definition:

---

[1] Gilbarg, David; Trudinger, Neil S. Elliptic partial differential equations of second order. Reprint of the 1998 edition. Classics in Mathematics. Springer-Verlag, Berlin, 2001.

Let function $f: D \rightarrow \{A_1, A_2,...,A_n\}$ and let $A_1 < A_2 < ... < A_n$. If $a$ and $b$ are adjacent in $D$, assume $f(a)=A_i$, then $f(b)= A_i$, $A_{i-1}$ or $A_{i+1}$. Point $(a,f(a))$ and $(b,f(b))$ are then said to be gradually varied.

A 2D function (surface) is said to be gradually varied if every adjacent pair is gradually varied.

Discrete Surface Fitting Definition:

Given $J \subseteq D$, and $f: J \rightarrow \{A_1,A_2,...A_n\}$, decide if there exists an F: $D \rightarrow \{A_1,A_2,...,A_n\}$ such that $F$ is gradually varied where $f(x)=F(x)$, and $x$ is in $J$.

An example using real numbers is:

Let D be a subset of real numbers $D \rightarrow \{1,2,.,.,.,n\}$. If $a$ and $b$ are adjacent in D such that $|f(a)-f(b)| \leq 1$ then point $(a,f(a))$ and $(b,(f(b))$ is said to be gradually varied.

A 2D function (a surface) is said to be gradually varied if every adjacent pair are gradually varied.

There are three theorems necessary for this development:

**Theorem 2.2** (Chen, 1989) [8][9] The necessary and sufficient conditions for the existence of a gradually varied extension $F$ are: for all $x,y$ in $J$, $d(x,y) \geq |i-j|$, $f(x)=A_i$ and $f(y)=A_j$, where $d$ is the distance between $x$ and $y$ in $D$.

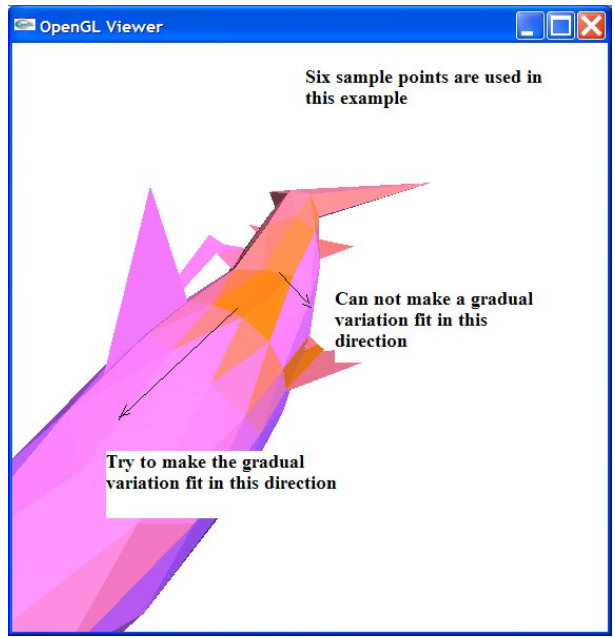

Fig. 2.1 The GVF algorithm using 6 guiding points (Original data from Princeton 3D Benchmark)

We have implemented the algorithms for smooth reconstruction in 2D rectangle domain and gradually varied reconstruction for 2D manifold [2][3]. In theory, McShane and Whitney obtained an important theorem for Lipschitz function extension [13][15]. Kirszbraun and later

Valentine studied the Lipschitz mapping extension for Hilbert spaces [16]. Recent theoretical research works on these aspects can be found in [11].

The gradually varied function is highly related to the Lipschitz function and the local Lipschitz function. It was proposed for discrete surface and data reconstruction [2,5].

## 3. New Algorithm Design and Analysis

Data reconstruction is used to fit a function based on the observations of some sample (guiding) points. The most popular application is to fit such a surface in a rectangle domain. This paper presents new developed algorithms to fit a function on a manifold (2D or 3D).

It is hard to make a smooth function on a closed manifold if do not allow manifold refinement. Some methods are invented in the past based on relatively dense sample points [1]. GVF will be only good for the continuous function; it is hard for the Taylor expansion in local with arbitrary neighborhood [10,12,14]. Local Euclidean block usually adopted by computer graphics in subdivision problem [12,14].

Our design is to partition the manifold into several components, so that each sample point is at the edge (boundary) of the partition. The edge of the component is a simple closed curve (path in terms of graph). Each closed curve contains original sample points. We also assume that the component is simply connected. Then find gradually varied interpolation on paths/curves, (each curve, one-by-one). So we will have data on all boundaries of all components; the fitted data points are all gradually varied. In each component, inside of each component does not contain any data or sample points except boundary edges. Finally, we use harmonic function to fit each component using a simple iterating process. So we will have a locally harmonic function on the surface. If we need a smooth function on the surface, we only need to apply a subdivision algorithm to the above solution.

### 3.1 The Main Algorithm Design

We will first present this algorithm in major steps, and then we will give detailed information for each step.

**Algorithm 3.1. Major Steps**
**Input:** Discretization of a 3D manifold (in 3D space) *M* and its boundary *D*. Some randomly arranged sample points are located on the boundary surface.
**Output**: Get a smooth extension based on sample points on the boundary surface; get the smooth extension for the whole 3D manifold.

**Step 1.** *D* is represented by a simplicial complex form (triangulated format) or cubical form. Sample points are located randomly in *D*. Find *A1, A2,...,Am* in Real Number *R*, so that the sample points on *D* satisfy the gradual variation condition.

**Step 2.** Partition *D* into simply connected components, so that each sample point is at the edge (boundary) of the partition.

**Step 3.** Extract each boundary curve for each component. Each boundary curve is a simple closed curve (path). If there is a pathological case, correct it using an interactive process. We get all sample points on the boundary curve.

**Step 4.** Get a gradually varied interpolation on boundary curves obtained in Step 3.

**Step 5.** Random assign values to each inner points of each component. Or use gradually varied fitting to get the initial value of the inner points in each component (a fast way to converge for the next step).

**Step 6.** Harmonic data reconstruction: Update the value in inner points based on the average value of neighboring points. This process is for each components until converge.

Add the smooth subdivision algorithm as needed.

**Step 7.** Use boundary surface calculated above for all inner points of 3D manifold M. Random assign values to each inner points of M.

**Step 8.** Harmonic data reconstruction for 3D: Update the value in inner points based on the average value of neighboring points. This process is for each components until converge.

Now, we explain the detailed implementation for each step. The main challenger is Step 2.

**Algorithm for Step 2**: There will be two sub-algorithms. **Sub-algorithm 1:** The simplest algorithm will be to link to all points one by one, until form a cycle. Since we will have finite number of points. The process always can find a cycle. If a simple cycle is found that already contain all sample points, we just use it for the data reconstruction. See the example in next section.

Other wise start at any point in the fixed cycle. Link to a sample points, until a cycle is formed or contains all the rest of sample points, then we link the last sample points to the existing cycle. This process will stop. This is a linear time algorithm using a queue.

**Sub-algorithm 2:** Geodesic-like partition. Based on two points near, find a large cycle. Select a point in the cycle and select a sample point to find a large cycle. Until all sample points are included. The process will stop.

**Lemma 3.1.** For any closed 2D simply connected discrete surfaces, if all points on the edge of a partition are connected, then each partitioned patch is simply connected.

*Proof*: This is because that if there is a hole in a partitioned patch, the edge of the hole will be isolated. Therefore, it is not connected to other edges of the partition.

**Lemma 3.2.** For any closed 2D discrete surfaces with genus $g$, if all sample points are connected by gradually varied paths, (each of vertex in the path has two neighbor vertices in the path) then the paths are the boundary of several manifolds.

Based on the Theorems provided in Section 2, we will have the harmonic extension for the partitioned manifolds.

## 3.2 The Algorithm Analysis

In this subsection, we will perform some rough algorithm analysis to the algorithms designed in above subsection.

**As we discussed the input of the main algorithm is:** (1) M the 3D manifold, (2) D the boundary of M, a closed two D surface in the discrete format. Let's assume D is triangulated or a set of polygon in 3D (2D cell complexes), (3) J the sample points which is the subset of D with assigned values in R.

There are two types of distances: Geodesic distance of two points on the boundary of M, and digital-discrete distance of the two points on D, i.e., the graph-distance, how many vertices in between from one vertex to another. Geodesic distance on the edge if there is an edge in D; (otherwise, we just use the shortest path to approximate the real geodesic distance.) By the way, we will refer the distance as graph-distance not geodesic-distance, otherwise will specify specifically.

**The output of the main algorithm**: a smooth extension based on sample points on the boundary surface; and then the smooth extension for the whole 3D manifold M.

**Lemma 3.3.** The time complexity for complete Step 1 in Algorithm 3.1 is $O(|J| |D| \log|D|)$.

*Proof:* Calculate distance for each pair of points in *J*. Just like to perform *|J|* times the Dijstra algorithm. The Dijstra algorithm uses $O(|E|+|V|\log |V|)$ time [20]. Where |E| and |V| are numbers of edges and vertices. Since $|E| = O(|V|)$ in D. So the lemma is proven.
We can get both distances types. Next, we calculate the "tangent" or "slop" of each pair if we just want to just get a Lipschitz extension for the initial values. For none Lipschitz extension, we need to use algorithms described in [2].

**Lemma 3.4.** The time complexity for complete Step 2-3 in Algorithm 3.1 is $O(|D|)$.

*Proof:* The Step 2 is to partition *D* into simply connected components, so that each sample point is at the edge (boundary) of the partition.

A cording to Step, we can use Geodesic distance to find a line from a point to another point in J. This way will eliminate some of the paths that are the same digital-discrete distance. The algorithm will need to find two difference lines from each point visited. This is a linear algorithm about |D| since it is possible that every points in D will be included. This is because we got the distance (weighted) and just fill the road to each, when a point is visited, we will mark it in a queue. And then continue to get the objective points. When a visited point has more than two neighbors are visited, we will delete this point from the queue.

This process can be combined with the process in Step 4.

**Lemma 3.5.** The time complexity for complete Step 4 in Algorithm 3.1 is $O(|D| |D|)$. And the improved algorithm can reach the $O(|D| \log |D|)$.

*Proof:* This step is to assign all values on the marked points using gradually varied filling algorithm. It usually takes $O(|D|^2)$ if a simple algorithm is implemented.

This algorithm can be implemented in $O(|D| \log|D|)$ since the decomposition of the mesh is Jordan. [7]

**Lemma 3.6.** The time complexity for complete Step 5-6 in Algorithm 3.1 is $O(|D|^3)$ if we use the Gaussian elimination method . It can usually reach $O(|D| \log |D|)$ if we use fast approximation algorithm for sparse symmetrically diagonally dominant linear systems .

*Proof:* Two kinds of implementations: iteration or solving the system of linear equations. (1) Random assign values to each inner points of each component. Or use gradually varied fitting to get the initial value of the inner points in each component (a fast way to converge for the next step).

(2) This stage can also be implemented by solving a system of linear equations. The linear equations can be formed in the following way:

All unknown vertices are named as $x_1,...x_N$. For each xi, we got its neighbors, some of them are in $x_1,...x_N$, and the rest of them are known (meaning that we know their values.) According to the definition of the equation of discrete harmonic functions:

$$a_{i1} x_1 + a_{i2} x_2 + ... + a_{ii} x_i + ... + a_{iN} x_N = C_i$$

where $a_{ii}$ equals to the integer indicating the number of neighbors. $a_{ik}$ is either equals to "-1" or "0." Where "-1" means that it is $x_k$ is a neighbor of $x_i$. $C_i$ is a real number that is the summation of values that are from neighbors not in the set of $x_1,...,x_N$.

In most of cases, the number of neighbors is bounded. And its matrix of coefficients $a_{ij}$ is sympatric since if $a_{ik}$, $k$ is not equal to $i$, is "-1" that means $x_k$ is a neighbor of $x_i$, so $x_i$ is a neighbor of $x_k$, that is to say $a_{ki}$ is also "-1."

If we use a cubic space in 3D, there will be Only 6-26 none-zeros in a row since a voxel only has at most 26 adjacent voxels.

This system of linear equations is a sparse symmetrically diagonally dominant linear system usually with a constant band-width in diagonally meaning that constant number of elements are not zero in a row. [20-21] The recent result by Koutis et al has reached the optimal to get the solution $O(m \log n \log(1/epsilon))$, where epsilon is a given constant for accuracy measurement [21].

**Lemma 3.7.** The time complexity for complete Step 7-8 in Algorithm 3.1 is $O(|D|^3)$ if use the Gaussian elimination method . It can usually reach $O(|D| \log |D|)$ if we use fast approximation algorithm for sparse symmetrically diagonally dominant linear systems.

*Proof:* The same as Lemma 3.8. We just implement the same algorithm to fill the inner points of a 3D with the values on the boundary surfaces. This algorithm does not need to assume the convexity of a geometric object.

Therefore,

**Theorem** 3.1. There is a near $O(|D|\log|D|)$ algorithm (Step 1-6 in Algorithm 3.1) for a Lipschitz extension to get a piecewise harmonic function on a 2D closed surface.

**Theorem** 3.2. There is a near O(|D|log|D|) algorithm for a three dimensional fill.

**3.3 Open Problems**

More sophisticated algorithms can be done in future research. We have three Open Problems are listed below:

**Open Problems:**
(1) The Balanced partition: each component are about the sample size. Is the time complexity in P or NP-hard?
(2) The minimum length partition: The length of all curves is minimum. This problem in cubical case might not be NP hard.
(3) The largest numbers of inner points. Similar to the second problem. More smooth part in fitted function. Is this problem NP-hard?

## 4. Implementation and Examples

The implementation is based on the basic property of discrete harmonic functions for the 2D closed manifold. A closed curve partition will be needed to get a unique function. As long as the function on the curve (path) are continuous (gradually varied). So the inside of the curve will be uniquely defined based on the harmonic function.

Data reconstruction is used to fit a function based on the observations of some sample (guiding) points.

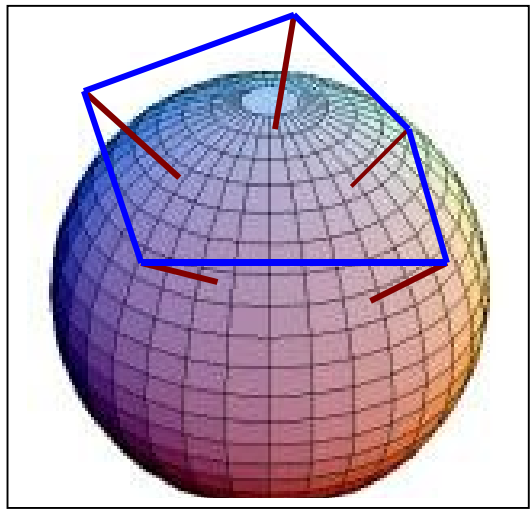

Fig. 4.1 An example of finding a closed path of cells for 2D manifold

For instance, we have five sample points on a sphere in Fig.1 . The height of line segment at a cell indicates the value of the sample point. So we can link all sample points to form a closed curve. Its corresponding path on the sphere separates the domain into two regions or components.

We have discovered a smooth data reconstruction method to reconstruct a function in 2D region based on the digital-discrete method [2,3]. The related method is based on gradually varied surface fitting [5-8]. This method differs from some successful methods have been discovered or proposed to solve the problem including the Voronoi-based surface method and the moving least square method [1][10][12][14].

We have implemented the harmonic fitting based on one closed boundary curve. See Fig. 4.2. The GVF method can help us to determine an initial value of the harmonic fitting. See Fig. 2

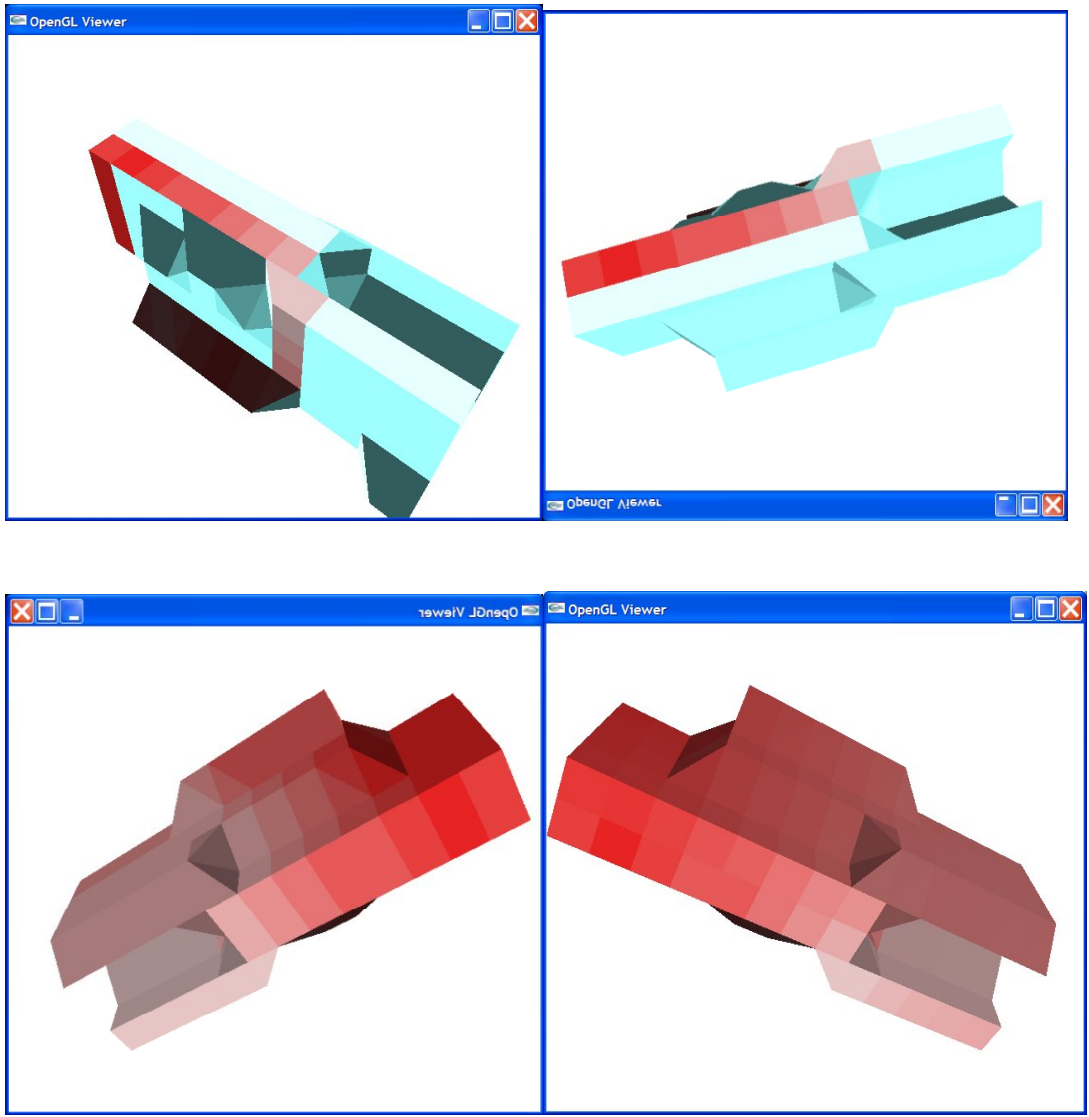

Fig. 4. 2. The selected cells form a boundary curve that is gradually varied: (a) and (b) are two displays for the guiding points (cells). (c) The GVF result. (d) The Harmonic fitting based on GVF (100 iterations).

The difference of gradually varied fitting and harmonic fitting was discussed in [9]. The following figure also show some similarity for directly fitting based on discrete sample points (not a path) .

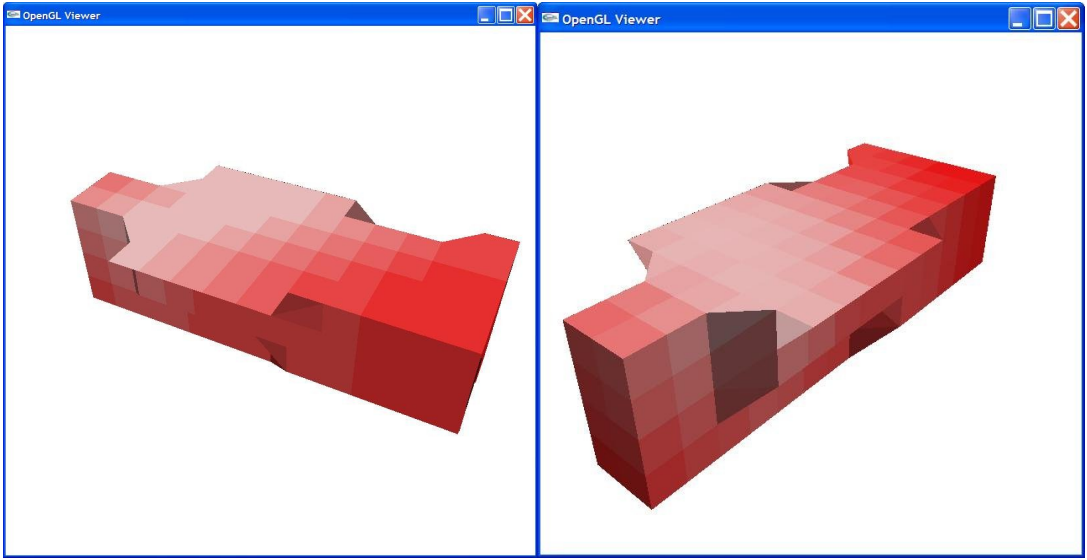

Fig. 4.3. Using 7 Points to fit the data on 3D surface: (a) the GVF result. (b) The Harmonic fitting (iterations few times) based on GVF.

When multiple cycles are constructed for the surface decomposition, we still can get a very good harmonic fitting.

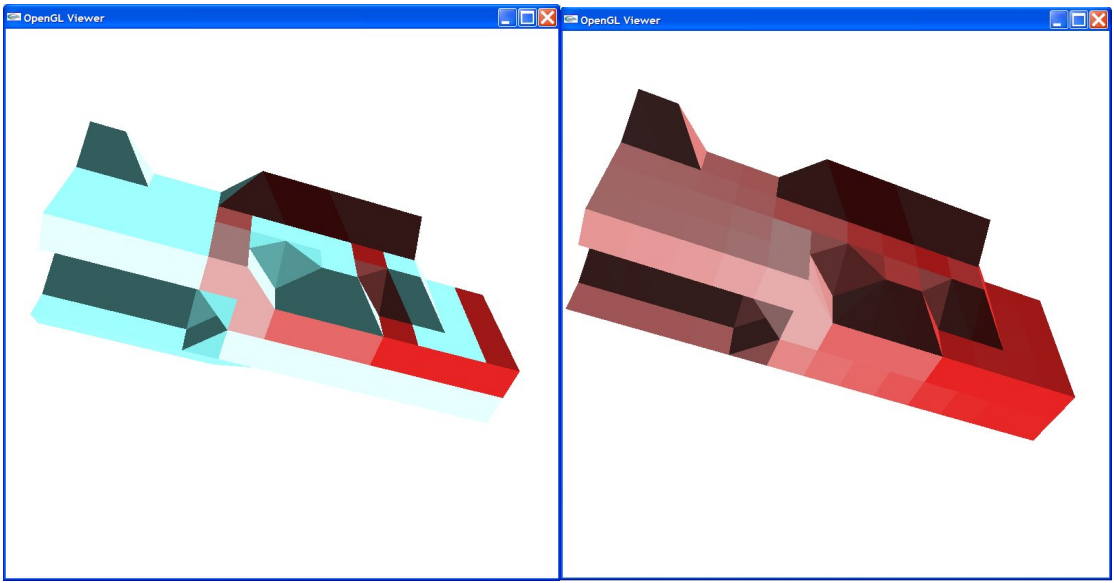

Fig. 4.4. Three components of the 2D manifolds: (a) the GVF on the paths. (b) The Harmonic fitting based on GVF as guiding points.

The following example shows the guiding points that are not connected initially. But our algorithm will find the paths to connect those points. We have 12 guiding points in this example. Then connecting all guiding points will result a petition. The harmonic reconstruction of the surface on surface based on the fitted paths will be shown in Fig. 4.7.

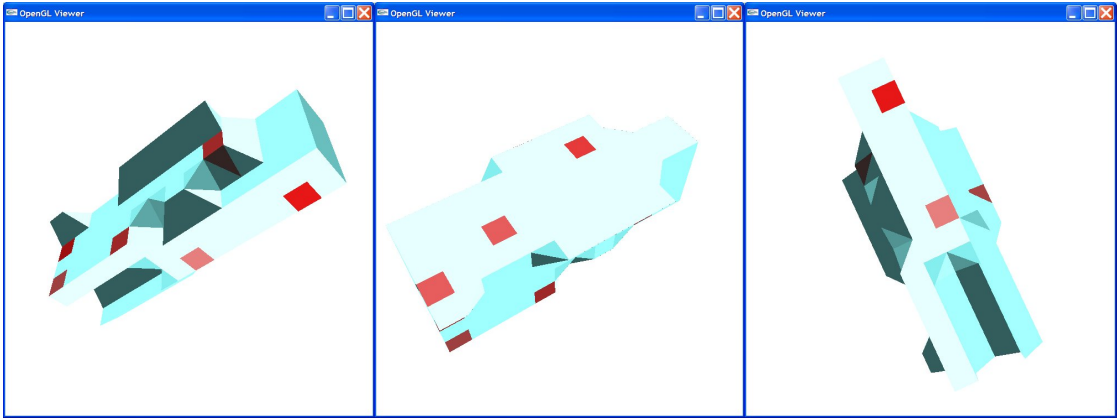

Fig. 4.5. Twelve guiding points from difference view angles.

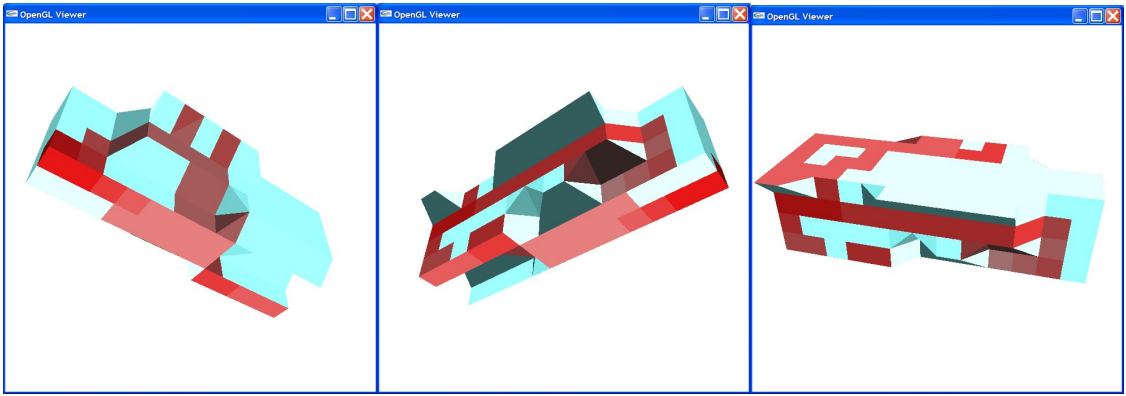

Fig. 4.6. The gradually varied paths for partition based on the twelve guiding points.

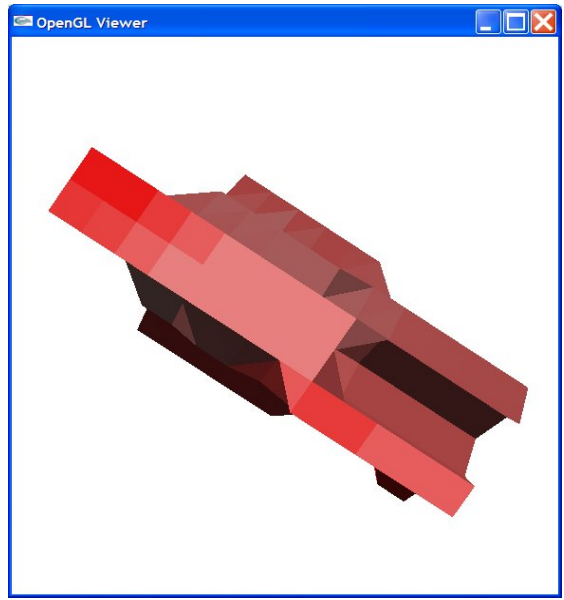

Fig. 4.7. The harmonic reconstruction of the surface on surface based on the reconstructed paths.

In real data practice, in order to make an easy implementation, we could first do the gradually varied surface fitting on the boundary surface. Then we can find all partitions. Finally we do the harmonic function calculation.

## 5. Discussions

The classical discrete method for data reconstruction is based on domain decomposition according to guiding (or sample) points. Then the Spline method (for polynomial) or finite elements method (for PDE) is used to fit the data. Our method is based on the gradually varied function that does not assume the property of being linearly separable among guiding points, i.e. no domain decomposition methods are needed. We also demonstrate the flexibility of the new method and its potential to solve a variety of problems.

Our approach is a new treatment to the important problem in the applications of numerical computation and computer graphics. There are some popular cited research works for data reconstruction. Levin used the moving least square (MLS) method for mesh-free surface interpolation [12][17]. This method still need to use dense sample points meaning that if you do not have few points in your neighborhood, you need to adjust the size of your neighbor. This is not an automated solution. An adaptive process is needed, which requires an artificial intelligence resolution. Our method tries to find the first fit for MLS or other subdivision methods.

Bajaj et al used Delaunay triangulations for automated reconstruction of three D objects [18]. Again, uses the linear interpolation as the first treatment.

Harmonic functions, on the other hand, are not limited to linear functions. But all linear functions are harmonic. Therefore, harmonic functions can use a polygon as the boundary, not only triangles and rectangles.

This paper mainly presents the results using harmonic functions on closed 2D manifolds. This method can be easily extended to higher multi-dimensions. We also include an advanced consideration related to the use of gradually varied mapping.

*Acknowledgements:* Professor Thomas Funkhouser at Princeton University provided helps on the 3D data sets and OpenGL display programs.